\documentclass[leqno,11pt]{amsart}
\usepackage{amsmath,amstext,amssymb,amsopn,amsthm,mathrsfs}

\theoremstyle{plain}
\newtheorem{teo}{Theorem}[section]
\newtheorem{defi}{Definition}[section]

\newtheorem{lemma}{Lemma}[section]
\newtheorem{proposition}{Proposition}[section]
\newtheorem{obs}{Observation}[section]
\numberwithin{equation}{section}
\newcommand{\dem}{\medskip \par \noindent \mbox{\bf Proof. }}

\def\ep{\hfill{$\Box $}}

\begin{document}

\title[Riesz \& Bessel Potentials and Fractional Derivatives ] {Riesz Potentials, Bessel Potentials and Fractional Derivatives  on Triebel-Lizorkin spaces for the Gaussian Measure.}

\author{A. Eduardo Gatto}
\address{Department of Mathematical Sciences, DePaul University, Chicago, IL,
   60614, USA.}
   \email{aegatto@depaul.edu}
   \author{Ebner Pineda}
\address{Departamento de Matem\'{a}tica,  Decanato de Ciencia y
Tecnologia, UCLA
 Apartado 400 Barquisimeto 3001 Venezuela.}
\email{epineda@ucla.edu.ve}
\author{Wilfredo~O.~Urbina}
\address{Department of Mathematics and Actuarial Sciences, Roosevelt University, Chicago, IL,
   60605, USA.}
\email{wurbinaromero@roosevelt.edu}



\begin{abstract}
In \cite{gatpiur} the  boundedness properties of Riesz Potentials, Bessel potentials and Fractional
Derivatives were studied in detail on Gaussian Besov-Lipschitz spaces $B_{p,q}^{\alpha}(\gamma_d)$. In this paper we will continue our study proving the boundedness of those operators on Gaussian  Triebel-Lizorkin spaces $F_{p,q}^{\alpha}(\gamma_d)$. Also these results can be extended to the case of Laguerre or Jacobi expansions and even further to the general framework of diffusions semigroups.
\end{abstract}

\maketitle

\section{Introduction and Preliminaries}

On $\mathbb{R}^d$ let us consider the Gaussian measure 
\begin{equation}
\gamma_d(dx)=\frac{e^{-\|x\|^2}%
}{\pi^{d/2}} dx, \, x\in\mathbb{R}^d
\end{equation}
 and the Ornstein-Uhlenbeck
differential operator
\begin{equation}\label{OUop}
L=\frac12\triangle_x-\left\langle x,\nabla _x\right\rangle.
\end{equation}

Let $\nu=(\nu _1,...,\nu_d)$ be a
multi-index such that $\nu _i \geq 0, i= 1, \cdots, d$,  let $\nu
!=\prod_{i=1}^d\nu _i!,$ $\left| \nu \right| =\sum_{i=1}^d\nu _i,$ $%
\partial _i=\frac \partial {\partial x_i},$ for each $1\leq i\leq d$ and $%
\partial ^\nu =\partial _1^{\nu _1}...\partial _d^{\nu _d},$ consider the normalized Hermite polynomials of order $\nu$ in $d$ variables,
\begin{equation}
h_\nu (x)=\frac 1{\left( 2^{\left| \nu \right| }\nu
!\right)
^{1/2}}\prod_{i=1}^d(-1)^{\nu _i}e^{x_i^2}\partial _i^{\nu _i}(e^{-x_i^2}),
\end{equation}
it is well known, that the Hermite polynomials are
eigenfunctions of the operator $L$,
\begin{equation}\label{eigen}
L h_{\nu}(x)=-\left|\nu \right|h_\nu(x).
\end{equation}
Given a function $f$ $\in L^1(\gamma _d)$ its
$\nu$-Fourier-Hermite coefficient is defined by
\[
\hat{f}(\nu) =<f, h_\nu>_{\gamma_d}
=\int_{\mathbb{R}^d}f(x)h_\nu (x)\gamma _d(dx).
\]
Let $C_n$ be the closed subspace of $L^2(\gamma_d)$ generated by
the linear combinations of $\left\{ h_\nu \ :\left| \nu
\right| =n\right\}$. By the orthogonality of the Hermite
polynomials with respect to $\gamma_d$ it is easy to see that
$\{C_n\}$ is an orthogonal decomposition of $L^2(\gamma_d)$,
$$ L^2(\gamma_d) = \bigoplus_{n=0}^{\infty} C_n,$$
this decomposition is called the Wiener chaos.

Let $J_n$ be the orthogonal projection  of $L^2(\gamma_d)$ onto
$C_n$, then if $f\in L^2(\gamma_d)$
\[
J_n f=\sum_{\left|\nu\right|=n}\hat{f}(\nu) h_\nu.
\]
Let us define the Ornstein-Uhlenbeck semigroup $\left\{
T_t\right\} _{t\geq 0}$ as
\begin{eqnarray}\label{01}
\nonumber T_t f(x)&=&\frac 1{\left( 1-e^{-2t}\right) ^{d/2}}\int_{\mathbb{R}^d}e^{-\frac{%
e^{-2t}(\left| x\right| ^2+\left| y\right| ^2)-2e^{-t}\left\langle
x,y\right\rangle }{1-e^{-2t}}}f(y)\gamma _d(dy)\\
& = & \frac{1}{\pi^{d/2}(1-e^{-2t})^{d/2}}\int_{\mathbb R^d} e^{-
\frac{|y-e^{-t}x|^2}{1-e^{-2t}}} f(y) dy
\end{eqnarray}
The family $\left\{ T_t\right\}_{t\geq 0}$ is a strongly
continuous Markov semigroup on $L^p(\gamma_d)$, $1 \leq p \leq
\infty$, with infinitesimal generator $L$. Also, by a change of
variable we can write,
\begin{equation}\label{t1}
T_t f(x)=\int_{\mathbb{R}^d} f(\sqrt{1-e^{-2t}}u + e^{-t}x)\gamma
_d(du).
\end{equation}

Now, by Bochner subordination formula, see Stein \cite{se70}, we
define the Poisson-Hermite semigroup $\left\{ P_t\right\} _{t\geq
0}$ as
\begin{equation}\label{PoissonH}
 P_t f(x)=\frac 1{\sqrt{\pi }}\int_0^{\infty} \frac{e^{-u}}{\sqrt{u}}T_{t^2/4u}f(x)du
\end{equation}
From (\ref{01}) we
obtain, after the change of variable $r=e^{-t^2/4u}$,
\begin{eqnarray}\label{03}
\nonumber P_t f(x)&=&\frac 1{2\pi
^{(d+1)/2}}\int_{\mathbb{R}^d}\int_0^1t\frac{\exp \left( t^2/4\log
r\right) }{(-\log r)^{3/2}}\frac{\exp \left( \frac{-\left|
y-rx\right| ^2}{1-r^2}\right) }{(1-r^2)^{d/2}}\frac{dr}rf(y)dy\\
&=& \int_{\mathbb{R}^d} p(t,x,y) f(y)dy,
\end{eqnarray}
with
\begin{equation}
p(t,x,y) = \frac 1{2\pi ^{(d+1)/2}}\int_0^1t\frac{\exp \left(
t^2/4\log r\right) }{(-\log r)^{3/2}}\frac{\exp \left(
\frac{-\left| y-rx\right| ^2}{1-r^2}\right)
}{(1-r^2)^{d/2}}\frac{dr}r.
\end{equation}
Also by the change of variable $s= t^2/4u$ we have,
\begin{equation}
P_t f(x)=\frac 1{\sqrt{\pi }}\int_0^{\infty} \frac{e^{-u}}{\sqrt{u}}T_{t^2/4u}f(x)du
=\int_0^{\infty} T_s f(x) \mu^{(1/2)}_t(ds),
\end{equation}
where the measure
\begin{equation}\label{onesided1/2}
\mu^{(1/2)}_t(ds) = \frac t{2\sqrt{\pi
}}\frac{e^{-t^2/4s}}{s^{3/2}}ds,
\end{equation}
is called the one-side stable measure on $(0, \infty)$ of order
$1/2$.

The family $\left\{ P_t\right\}_{t\geq 0}$ is a strongly
continuous semigroup on $L^p(\gamma_d)$, $1 \leq p \leq \infty$,
with infinitesimal generator $-(-L)^{1/2}$, furthermore $\left\{ P_t\right\}$ is an analytic semigroup in $t$. In what follows, often we are  going to use the notation $$u(x,t) = P_{t}f(x) \; \mbox{and} \; u^{(k)}(x,t) = \frac{\partial^{k}}{\partial t^k} P_{t}f(x).$$

Observe that by (\ref{eigen}) we have that
\begin{equation}\label{OUHerm}
T_t h_\nu(x)=e^{-t\left| \nu\right|}h_\nu(x),
\end{equation}
and
\begin{equation}\label{PHHerm}
 P_t h_\nu(x)=e^{-t\sqrt{\left| \nu\right|}}h_\nu(x),
\end{equation}
i.e. the Hermite polynomials are eigenfunctions of $T_t$ and $P_t$ for any $t \geq 0$.\\

The operators that we are going to consider in this paper are the following:
\begin{itemize}

\item  For $\beta>0$, the Fractional Integral or Riesz potential of order $\beta$, $I_\beta$, with respect to the Gaussian measure is defined formally as
\begin{equation}\label{i1}
I_\beta=(-L)^{-\beta/2}\Pi_{0},
\end{equation}
where,
$\Pi_{0}f=f-\displaystyle\int_{\mathbb{R}^{d}}f(y)\gamma_{d}(dy)$,
for $f\in L^{2}(\gamma_{d})$. That means that for  the Hermite
polynomials $\{h_\nu\}$, for $\left|\nu\right|>0$,
\begin{equation}\label{e4}
I_\beta h_\nu(x)=\frac 1{\left|
\nu \right|^{\beta/2}}h_\nu(x),
\end{equation}
and for $\nu=0=(0,...,0), \, I_{\beta}(h_{0})=0.$ Then by linearity, $I_\beta^{\gamma}$ can be defined to any polynomial.

In this paper the following representation of $I_\beta^{\gamma}$ is crucial to get the results. If $f$ is a polynomial,
\begin{equation}\label{e3}
I_\beta f(x)  =\frac 1{\Gamma(\beta)}\int_0^{\infty}
t^{\beta-1}(P_t f(x) -P_{\infty} f(x))\,dt,
\end{equation}
where $P_{\infty} f(x) = \int_{\mathbb{R}^{d}} f(y) \gamma_{d}(dy)$. The proof is very simple but we include it because we will refer to this calculation several times. By linearity it is enough to do it for the Hermite polynomials. Let $|\nu| >0$, otherwise the proof is trivial, then by orthogonality and the change of variable $u =t\sqrt{\left| \nu\right|}$, we have

$$\frac 1{\Gamma(\beta)} \int_{\mathbb{R}^{d}}  t^{\beta-1}(P_t h_\nu(x) -P_{\infty} h_\nu(x)) dt \hspace{3cm}$$
\begin{eqnarray*}
&& \hspace{1cm} =   \frac 1{\Gamma(\beta)} \int_{\mathbb{R}^{d}}  t^{\beta-1} e^{-t\sqrt{\left| \nu\right|}}dt h_\nu(x)\\
&& \hspace{1cm} =     \frac 1{\Gamma(\beta)} \int_{\mathbb{R}^{d}} (\frac{ u}{\sqrt{\left| \nu\right|}})^{\beta-1} e^{-u}\frac{du}{\sqrt{\left| \nu\right|}} h_\nu(x)= \frac 1{\left| \nu \right|^{\beta/2}} \, h_\nu(x).
\end{eqnarray*}

\begin{obs}\label{LpRiesz}  P. A. Meyer's multiplier theorem shows that $I_\beta$ has a continuous extension  to $ L^p(\gamma _d)$,  $1 < p < \infty$, and then (\ref{e3}) can be extended for $f \in L^{p}(\gamma_d)$, by the density of the polynomials, see \cite{berg}.
\end{obs}

\item  The Bessel Potential of order $\beta>0,$
$\mathcal{J}_\beta$, associated to the Gaussian  measure
is defined formally as
\begin{eqnarray}
\mathcal{J}_\beta= (I+\sqrt{-L})^{-\beta},
\end{eqnarray}
meaning that for the Hermite polynomials we have,
\begin{eqnarray*}
\mathcal{J}_\beta h_\nu(x)=\frac 1{(1+\sqrt{\left|
\nu\right|})^{\beta}}h_\nu(x).
\end{eqnarray*}
Again  by linearity $\mathcal{J}_\beta$ can be extended to any polynomial. By similar calculation as above (\ref{e3}),  the Bessel potentials can be represented 
\begin{equation}\label{Besselrepre}
\mathcal{J}_\beta
f(x)=\frac{1}{\Gamma(\beta)}\int_{0}^{+\infty}t^{\beta}e^{-t}P_{t}f(x)\frac{dt}{t},
\end{equation}

\begin{obs}\label{LpBessel}   P. A. Meyer's theorem shows that $\mathcal{J}_\beta$  is a \mbox{continuous} operator on $L^p(\gamma_d),$ $1 < p < \infty$, and again (\ref{Besselrepre}) can be extended for $f \in L^{p}(\gamma_d)$, by the density of the polynomials there.
\end{obs}

\item The Riesz fractional derivate of order $\beta >0$ with respect to the
Gaussian measure $D^\beta$, is defined  formally as
\begin{equation}
D^\beta=(-L)^{\beta/2},
\end{equation}
meaning that for the Hermite polynomials, we have
\begin{equation}\label{e6}
D^\beta h_\nu(x)=\left| \nu\right|^{\beta/2}
h_\nu(x),
\end{equation}
 thus by linearity can be extended to any polynomial, see  \cite{lu} and \cite{ebner}.
 
 Now, if $f$ is a polynomial, by the linearity of the operators
$I_\beta$ and $D^\beta$, (\ref{e4}) and
(\ref{e6}), we get
\begin{equation}\label{di}
 \Pi_0f=I_\beta(D^\beta f)=D^\beta(I_\beta f).\\
\end{equation}

In the case of $0 < \beta <1$ we have the following integral representation,  for $f$ a polynomial,
\begin{equation}\label{e5}
D^\beta f(x) =\frac 1{c_\beta}\int_0^{\infty}t^{-\beta-1}(P_t  -I)\, f(x)  dt,
\end{equation} 
where
$
c_\beta=\int_0^\infty u^{-\beta-1}(e^{-u}-1)du.$

 Now, if $\beta \geq 1$, let  $k$   be the smallest integer greater than $\beta$ i.e. $ k-1 \leq \beta < k$, then the fractional derivative $D^\beta$ can be represented as 
\begin{equation}
D^\beta f = \frac{1}{c^k_{\beta}}\int_0^{\infty} t^{-\beta-1} ( P_t -I )^k f \, dt,
\end{equation}
where $c^k_{\beta} = \int_0^{\infty} u^{-\beta-1} ( e^{-u}-1 )^k  \, du$ and $f$ a polynomial, see \cite{sam}.

\item We also define a Bessel fractional derivative ${\mathcal D}^\beta$,  defined formally as
\[
{\mathcal D}^\beta=(I+\sqrt{-L})^{\beta},
\]
which means that for the Hermite polynomials, we have
\begin{equation}\label{e7}
{\mathcal D}^\beta h_\nu(x)=(1+ \sqrt{\left| \nu\right|})^{\beta} h_\nu(x),
\end{equation}
In the case of $0 < \beta <1$ we have the following integral representation,  
\begin{equation}\label{e8}
{\mathcal D}^\beta f =\frac 1{c_\beta}\int_0^{\infty}t^{-\beta-1}( e^{-t} P_t -I) \, f dt, 
\end{equation} 
where, as before, 
$
c_\beta=\int_0^\infty u^{-\beta-1}(e^{-u}-1)du$ and $f$ a polynomial.

Moreover,  if $\beta \geq1$ let $k$  be the smallest integer greater than $\beta$ i.e. $ k-1 \leq \beta < k$, we have the following representation of  ${\mathcal D}^\beta f$ 
\begin{equation}
{\mathcal D}^\beta f =  \frac{1}{c^k_{\beta}} \int_0^{\infty} t^{-\beta-1} (e^{-t} P_t -I )^k\,  f \, dt,
\end{equation}
where $c^k_{\beta} =  \int_0^{\infty} u^{-\beta-1} (e^{-u} -1 )^k du$ and $f$ a polynomial, see \cite{sam}.\\

\end{itemize}

The Gaussian Triebel-Lizorkin $F_{p,q}^{\alpha}(\gamma_d)$
spaces were introduced in \cite{piur02}, see also \cite{ebner}, as follows

 \begin{defi}
Let $\alpha \geq  0$, $k$ be the smallest integer
greater than $\alpha$, and $1\leq p,q<\infty$.  The Gaussian  Triebel-Lizorkin
 space $F_{p,q}^{\alpha}(\gamma_{d})$  is the set of functions $f\in
L^{p}(\gamma_{d})$ for which
\begin{equation}\label{e16}
\left\| \left( \int_0^{\infty} (t^{k-\alpha}
\left|\frac{\partial^{k}P_t f}{\partial t^{k}} \right|) ^{q}
\frac{dt}{t} \right) ^{1/q} \right\|_{p,\gamma_d}<\infty.
\end{equation}
The norm of $f \in F_{p,q}^{\alpha}(\gamma_d)$ is defined as
\begin{equation}
\left\| f \right\|_{F_{p,q}^{\alpha}}: =  \left\| f \right\|_{p,
\gamma_d} +\left\| \left( \int_0^{\infty} (t^{k-\alpha}
\left|\frac{\partial^{k}P_t f}{\partial t^{k}} \right|) ^{q}
\frac{dt}{t} \right) ^{1/q} \right\|_{p,\gamma_d}.
\end{equation}
\end{defi}
The definition of $F_{p,q}^{\alpha}(\gamma_d)$ does not depend on which $k>\alpha$ is chosen and the resulting norms are equivalent,  for the proof of this result and other properties of these spaces see \cite{piur02}.\\

For Gaussian Triebel-Lizorkin spaces we have the following inclusion result,

\begin{proposition}\label{incluTriebel}
The inclusion $F_{p,q_1}^{\alpha_{1}}(\gamma_d)\subset
F_{p,q_2}^{\alpha_{2}}(\gamma_d)$ holds for $\alpha_1 > \alpha_2>0$ and
$q_1 \geq q_2$.
\end{proposition}
For the proof of this see \cite{piur02}, Proposition 2.2.\\

In what follows we will need Hardy's inequalities, so for completeness we will write then here, see \cite{se70} page 272,
\begin{equation}\label{hardy1}
\int_{0}^{+\infty}\big(\int_{0}^{x}f(y)dy\big)^{p} x^{-r-1} dx \leq \frac{p}{r}\int_{0}^{+\infty}(y f(y))^{p}y^{-r-1}dy,
\end{equation}
and
\begin{equation}\label{hardy2}
\int_{0}^{+\infty}\big(\int_{x}^{\infty}f(y)dy\big)^{p} x^{r-1}dx \leq \frac{p}{r}\int_{0}^{+\infty}(y f(y))^{p}y^{r-1}dy,
\end{equation}
 where $f\geq 0, p\geq 1$ and $r>0.$\\

 Finally, in \cite{gatur} Gaussian Lipchitz spaces $Lip_{\alpha}(\gamma_d)$ were considered and the boundedness of  Riesz Potentials, Bessel potentials and Fractional
Derivatives on them and  in \cite{gatpiur} the boundedness of those operators were studied on Gaussian Besov-Lipschitz spaces $B_{p,q}^{\alpha}(\gamma_{d})$. In the next section we are going to
 study the boundedness properties of those operators for Gaussian Triebel-Lizorkin spaces $F_{p,q}^{\alpha}(\gamma_{d})$.
As usual in what follows $C$ represents a constant that is not necessarily the same in each occurrence.


\section{Main results and Proofs}
\medskip

 \begin{teo}\label{PotRiesz}
Let $\alpha\geq 0, \beta>0$, $1<p<\infty, 1\leq q < \infty$ then $I_{\beta}$ is
bounded from $F_{p,q}^{\alpha}(\gamma_{d})$ into
$F_{p,q}^{\alpha+\beta}(\gamma_{d})$.
\end{teo}
\dem 

Let $k > \alpha+\beta +1$ an integer fixed and $f \in F_{p,q}^{\alpha}(\gamma_{d}).$
Using the integral representation of Riesz Potentials (\ref{e3}), the semigroup property of $\{P_{t}\}$ and the fact that $P_{\infty}f(x)$ is a constant, we get
\begin{eqnarray}\label{est1}
\nonumber P_{t}(I_{\beta}f)(x)&=&\displaystyle\frac{1}{\Gamma(\beta)}\int_{0}^{+\infty}s^{\beta-1}P_{t}(P_{s}f(x)-P_{\infty}f(x))ds\\
&=&\displaystyle\frac{1}{\Gamma(\beta)}\int_{0}^{+\infty}s^{\beta-1}(P_{t+s}f(x)-P_{\infty}f(x))ds.
\end{eqnarray}
Then using again that $P_{\infty}f(x)$ is a constant and the chain rule,
\begin{eqnarray}\label{est2}
\nonumber \frac{\partial^{k}}{\partial t^{k}}P_{t}(I_{\beta}f)(x) &=&\frac{1}{\Gamma(\beta)}\int_{0}^{+\infty}s^{\beta-1}\frac{\partial^{k} }{\partial t^{k}}(P_{t+s}f(x)-P_{\infty}f(x))ds\\
&=& \frac{1}{\Gamma(\beta)}\int_{0}^{+\infty}s^{\beta-1}u^{(k)}(x, t+s) ds.
\end{eqnarray}

\begin{enumerate}

\item[i)] Case $\beta\geq 1$: Using (\ref{est2}), the change of variable $r=t+s$, $dr=ds$ and Hardy's  inequality (\ref{hardy2}), we have
\begin{eqnarray*}
&&\big(\int_{0}^{+\infty}\big(t^{k-(\alpha+\beta)}|\frac{\partial^{k}P_{t}(I_{\beta}f)(x)}{\partial
t^{k}}|\big)^{q}\frac{dt}{t}\big)^{1/q} \quad\quad \quad \quad \quad
\quad\quad \quad\quad \quad\\
&&\quad \quad \quad \quad \quad \leq \frac{1}{\Gamma(\beta)}\big(\int_{0}^{+\infty}t^{q(k-(\alpha+\beta))}\big(\int_{0}^{+\infty}s^{\beta-1} |u^{(k)}(x, t+s) | ds \big)^{q}\frac{dt}{t}\big)^{\frac{1}{q}}\\
&&\quad \quad \quad \quad \quad = \frac{1}{\Gamma(\beta)}\big(\int_{0}^{+\infty}t^{q(k-(\alpha+\beta))}\big(\int_{t}^{+\infty}(r-t)^{\beta-1}
|u^{(k)}(x,r)|dr\big)^{q}\frac{dt}{t}\big)^{\frac{1}{q}}\\
&&\quad \quad \quad \quad \quad \leq\frac{1}{\Gamma(\beta)}\big(\int_{0}^{+\infty}t^{q(k-(\alpha+\beta))} \big(\int_{t}^{+\infty}r^{\beta-1}
|u^{(k)}(x,r)|dr\big)^{q} \frac{dt}{t}\big)^{\frac{1}{q}}\\
&&\quad \quad \quad \quad \quad \leq \frac{1}{\Gamma(\beta)}\frac{1}{k-(\alpha+\beta)}
\big(\int_{0}^{+\infty}\big(r^{k-\alpha}
|u^{(k)}(x,r)|\big)^{q}\frac{dr}{r}\big)^{\frac{1}{q}},
\end{eqnarray*}
and therefore
\begin{eqnarray*}
&&\|\big(\int_{0}^{+\infty}\big(t^{k-(\alpha+\beta)}|\frac{\partial^{k}P_{t}(I_{\beta}f)}{\partial
t^{k}}|\big)^{q}\frac{dt}{t}\big)^{\frac{1}{q}}\|_{p,\gamma}\quad \quad\quad \quad \quad \quad
\quad\quad \quad\quad \quad\\
&&\quad \quad \quad \quad \quad \leq
C_{k,\alpha,\beta}
\|\big(\int_{0}^{+\infty}\big(r^{k-\alpha}
|\frac{\partial^{k}P_{r}f}{\partial
r^{k}}|\big)^{q}\frac{dr}{r}\big)^{\frac{1}{q}}\|_{p,\gamma} <\infty  ,
\end{eqnarray*}
since $f\in  F_{p,q}^{\alpha}$. By Observation \ref{LpRiesz} and the previous estimate, 
 $$\| I_{\beta}f\|_{F_{p,q}^{\alpha+\beta}} \leq C \| f\|_{F_{p,q}^{\alpha}}. $$ \\

\item[ii)] Case $0<\beta< 1$: Again using (\ref{est2}),
 \begin{eqnarray*}
&&\big(\int_{0}^{+\infty}\big(t^{k-(\alpha+\beta)}|\frac{\partial^{k}P_{t} (I_{\beta}f)(x)}{\partial
t^{k}}|\big)^{q}\frac{dt}{t}\big)^{\frac{1}{q}} \quad \quad\quad \quad \quad \quad
\quad\quad \quad\quad \quad\\
&&\quad \quad \quad \quad \quad \leq\frac{1}{\Gamma(\beta)}\big(\int_{0}^{+\infty}t^{q(k-(\alpha+\beta))}\big(\int_{0}^{+\infty}s^{\beta}
|u^{(k)}(x, t+s)|\frac{ds}{s}\big)^{q}\frac{dt}{t}\big)^{\frac{1}{q}} \\
&&\quad \quad \quad \quad \quad \leq
\frac{C}{\Gamma(\beta)}(\int_{0}^{+\infty}t^{q(k-(\alpha+\beta))-1}\big(\int_{0}^{t}s^{\beta-1}
|u^{(k)}(x, t+s)|ds\big)^{q}dt\big)^{\frac{1}{q}}\\
&&\quad \quad \quad \quad \quad \quad \quad+\frac{C}{\Gamma(\beta)}
(\int_{0}^{+\infty}t^{q(k-(\alpha+\beta))-1}\big(\int_{t}^{+\infty}s^{\beta-1}
|u^{(k)}(x, t+s)|ds\big)^{q}dt\big)^{\frac{1}{q}}\\
&&\quad \quad \quad \quad \quad=(I) + (II)  .
\end{eqnarray*}

Writing $s^{\beta-1} = s^{\frac{\beta -1}{q}+ \frac{\beta -1}{q'}},\, \frac{1}{q}+\frac{1}{q'} =1,$ and using H\"older's inequality in the internal integral,
\begin{eqnarray*}
(I) &\leq&\frac{C_{\beta}}{\beta^{q-1}} \big(\int_{0}^{+\infty}t^{q(k-\alpha)-\beta-1}\int_{0}^{t}s^{\beta-1}
|u^{(k)}(x, t+s)|^{q}ds\, dt\big)^{1/q}.
\end{eqnarray*}
By Fubini- Tonelli's theorem and using that $q(k-\alpha)-\beta-1>0$, as $k > \alpha + \beta +1$, we get
\begin{eqnarray*}
(I) &\leq& \frac{C_{\beta}}{\beta^{q-1}}\big(\int_{0}^{+\infty}s^{\beta-1}\int_{s}^{+\infty}t^{q(k-\alpha)-\beta-1}
|u^{(k)}(x, t+s)|^{q}dt \,ds\big)^{1/q}\\
&\leq& \frac{C_{\beta}}{\beta^{q-1}}\big(\int_{0}^{+\infty}s^{\beta-1} \int_{s}^{+\infty}(t+s)^{q(k-\alpha)-\beta-1}
|u^{(k)}(x, t+s)|^{q}dt \,ds\big)^{1/q}.
\end{eqnarray*}
Then, by
the change of variable  $r=t+s$ and Hardy's inequality (\ref{hardy2}),
\begin{eqnarray*}
(I)&\leq&\frac{C_{\beta}}{\beta^{q-1}} \big(\int_{0}^{+\infty}s^{\beta-1} \int_{2s}^{+\infty}r^{q(k-\alpha)-\beta-1}
|u^{(k)}(x, r)|^{q}dr \, ds\big)^{1/q}\\
&\leq&\frac{C_{\beta}}{\beta^{q-1}}\big( \int_{0}^{+\infty}s^{\beta-1} \int_{s}^{+\infty}r^{q(k-\alpha)-\beta-1}
|u^{(k)}(x, r)|^{q}dr \, ds \big)^{1/q}\\
&\leq&\frac{C_{\beta}}{\beta^{q-1}} \frac{1}{\beta^{1/q}} \big(\int_{0}^{+\infty}\big(r^{k-\alpha}|u^{(k)}(x, r)|\big)^{q}\frac{dr}{r}\big)^{1/q} .
\end{eqnarray*}\\
On the other hand, Since $\beta<1$ then $t<s$ implies that
$s^{\beta-1}<t^{\beta-1}$ and by the change of variable $r=t+s$ and Hardy's
inequality (\ref{hardy2}), as $k>\alpha+\beta+1>\alpha+1$, we obtain
\begin{eqnarray*}
  (II)  &\leq& C_{\beta} \big(\int_{0}^{+\infty}t^{q(k-\alpha-1)-1}\big(\int_{t}^{+\infty}
|u^{(k)}(x, t+s)|ds\big)^{q}dt\big)^{\frac{1}{q}}\\
&\leq&C_{\beta} \big(\int_{0}^{+\infty}t^{q(k-\alpha-1)-1}\big(\int_{2t}^{+\infty}
|u^{(k)}(x, r)|dr\big)^{q}dt\big)^{\frac{1}{q}}\\
&\leq&C_{\beta} \frac{1}{k-\alpha-1}\big(\int_{0}^{+\infty}\big(r^{k-\alpha}|\frac{\partial^{k}P_{r}f(x)}{\partial
r^{k}}|\big)^{q}\frac{dr}{r}\big)^{\frac{1}{q}}.
\end{eqnarray*}
Therefore,
\begin{eqnarray*}
&& \|\big(\int_{0}^{+\infty}\big(t^{k-(\alpha+\beta)}|\frac{\partial^{k}P_{t}(I_{\beta}f)(x)}{\partial
t^{k}}|\big)^{q}\frac{dt}{t}\big)^{\frac{1}{q}}\|_{p,\gamma}\\
& & \quad  \quad  \quad  \leq
C_{k,\alpha,\beta }\|\big(\int_{0}^{+\infty}\big(r^{k-\alpha}|\frac{\partial^{k}P_{r}f}{\partial
r^{k}}|\big)^{q}\frac{dr}{r} \big)^{\frac{1}{q}}\|_{p,\gamma} < \infty  .
\end{eqnarray*}
as $f\in  F_{p,q}^{\alpha}$. By Observation \ref{LpRiesz} and the previous estimate, we get
 $$\| I_{\beta}f\|_{F_{p,q}^{\alpha+\beta}} \leq C \| f\|_{F_{p,q}^{\alpha}}. $$  
 \ep
\end{enumerate}


Next we want to study the boundedness properties of the Bessel potentials on Triebel-Lizorkin spaces.
In \cite{piur02}, Theorem 2.4, the following result was proved, for completeness the proof will be given here too.
\begin{teo}\label{PotBessel}
Let $\alpha\geq 0, \beta>0$ then for $1< p < \infty$,    $1\leq q < \infty$ $\mathcal{J}_{\beta}$ is
bounded from $F_{p,q}^{\alpha}(\gamma_{d})$ into
$F_{p,q}^{\alpha+\beta}(\gamma_{d})$.
\end{teo}

\dem 

Let  $k>\alpha+\beta+1$  a fixed integer and  $f\in F_{p,q}^{\alpha}(\gamma_{d})$. Observe that
\begin{eqnarray*}
|\frac{\partial^{k}P_{t}}{\partial t^{k}}( {\mathcal{J}}_{\beta}f)(x)| &\leq& \frac{1}{\Gamma(\beta)}\int_{0}^{+\infty}s^{\beta}e^{-s}
|u^{(k)}(x,t+s)|\frac{ds}{s}\\
&\leq&\frac{1}{\Gamma(\beta)}\int_{0}^{+\infty}s^{\beta} |u^{(k)}(x,t+s)|\frac{ds}{s},
\end{eqnarray*}
and therefore
\begin{eqnarray*}
&&\big(\int_{0}^{+\infty}\big(t^{k-(\alpha+\beta)}|\frac{\partial^{k}}{\partial
t^{k}} P_{t}( {\mathcal{J}}_{\beta}f)(x)|\big)^{q}\frac{dt}{t}\big)^{1/q} \quad\quad \quad \quad \quad
\quad\quad \quad\quad \quad\\
&&\quad \quad \quad \quad \quad \leq \frac{1}{\Gamma(\beta)}\big(\int_{0}^{+\infty}t^{q(k-(\alpha+\beta))}\big(\int_{0}^{+\infty}s^{\beta}
|u^{(k)}(x,t+s)|\frac{ds}{s}\big)^{q}\frac{dt}{t}\big)^{\frac{1}{q}}.
\end{eqnarray*}
But this is the same term studied in the proof of Theorem \ref{PotRiesz}, so the argument (divided in two cases
$\beta\geq 1$ and $0<\beta< 1$) is totally analogous in this case. \ep\\


We will study next the boundedness of the Riesz fractional derivative $D^\beta$ on Triebel-Lizorkin spaces. 

\begin{teo}\label{DerRiesz<1}
Let $0<\beta<\alpha<1$, Let $1\leq p,q<\infty$ then
 $D^{\beta}$ is
bounded from $F_{p,q}^{\alpha}(\gamma_{d})$ into
$F_{p,q}^{\alpha-\beta}(\gamma_{d})$.
\end{teo}

\dem

 Let $f\in F_{p,q}^{\alpha}(\gamma_{d})$, using Hardy's inequality (\ref{hardy1}) with $p=1$, and the Fundamental Theorem of Calculus,
 \begin{eqnarray}\label{est4}
\nonumber |D^{\beta}f(x)| &\leq&\displaystyle\frac{1}{c_{\beta}}\int_{0}^{+\infty}s^{-\beta-1}|P_{s}f(x)-f(x)|ds\\
&\leq&
\nonumber \displaystyle\frac{1}{c_{\beta}}\int_{0}^{+\infty}s^{-\beta-1}\int_{0}^{s}|\frac{\partial
}{\partial r}P_{r}f(x)|dr \,ds\\
&\leq&\displaystyle\frac{1}{c_{\beta}\beta} \int_{0}^{+\infty}r^{1-\beta}|\frac{\partial
}{\partial r}P_{r}f(x)|\frac{dr}{r}  .
\end{eqnarray}
Thus, 
\begin{eqnarray}\label{est5}
\nonumber \|D^{\beta}f \|_{p,\gamma}
&\leq&\displaystyle C_{\beta} \|\int_{0}^{+\infty}r^{1-\beta}|\frac{\partial}{\partial r}P_{r}f| \frac{dr}{r}\|_{p,\gamma}\\&\leq& C_{\beta}  \|f\|_{F_{p,q}^{\alpha}}<\infty ,
\end{eqnarray}
because $ F_{p,q}^{\alpha}(\gamma_{d})\subset
F_{p,1}^{\beta}(\gamma_{d}) $ ($\alpha> \beta$ and $q\geq1$). Now, by analogous argument using Hardy's inequality (\ref{hardy1}) with $p=1$,
\begin{eqnarray*}
|\frac{\partial}{\partial t}P_{t}(D^{\beta}f)(x)|&\leq& \displaystyle\frac{1}{c_{\beta}}\int_{0}^{+\infty}s^{-\beta-1}|\frac{\partial}{\partial t}P_{t+s}f(x)- \frac{\partial}{\partial t}P_t f(x)| ds\\
&\leq& \displaystyle\frac{1}{c_{\beta}}\int_{0}^{+\infty}s^{-\beta-1} \int_0^{s} |u^{(2)}(x,t+r)| dr \,ds\\
&\leq& \displaystyle\frac{1}{c_{\beta}\beta }\int_{0}^{+\infty}r^{-\beta}  |u^{(2)}(x,t+r)| dr \,ds\\
\end{eqnarray*}

This implies that
\begin{eqnarray}\label{est10}
\nonumber \int_{0}^{\infty}\big(t^{1-(\alpha-\beta)}|\frac{\partial }{\partial
t}P_{t}(D_{\beta}f)(x)|\big)^{q}\frac{dt}{t}&\leq&\displaystyle\frac{1}{c_{\beta}\beta}
\int_{0}^{\infty}\big(t^{1-(\alpha-\beta)}\int_{0}^{+\infty}r^{-\beta}|u^{(2)}(x,t+r)|dr\big)^{q}\frac{dt}{t}\\
\nonumber &=&C_\beta \int_{0}^{\infty}\big(t^{1-(\alpha-\beta)}\int_{0}^{t}r^{-\beta}|u^{(2)}(x,t+r)|dr\big)^{q}\frac{dt}{t}\\
 && \quad + \,C_\beta \int_{0}^{\infty}\big(t^{1-(\alpha-\beta)}\int_{t}^{+\infty}r^{-\beta}|u^{(2)}(x,t+r)|dr\big)^{q}\frac{dt}{t}\\
\nonumber  &=&(I)+(II).
\end{eqnarray}

Writing $s^{\beta-1} = s^{\frac{\beta -1}{q}+ \frac{\beta -1}{q'}},\, \frac{1}{q}+\frac{1}{q'} =1,$ and using H\"older's inequality in the internal integral, we have
$$
(I)\leq C_\beta(1-\beta)^{1-q}\int_{0}^{\infty}t^{(2-\alpha)q-2+\beta} \int_{0}^{t}r^{-\beta}|u^{(2)}(x,t+r)|^{q}dr\,dt,
$$
then by Fubini-Tonelli's theorem, we get
 $$(I)\leq C_\beta(1-\beta)^{1-q}\int_{0}^{\infty}r^{-\beta} \int_{r}^{\infty}t^{(2-\alpha)q+\beta-2}|u^{(2)}(x,t+r)|^{q}dt\,dr .
$$
It is easy to prove that $(2-\alpha)q+\beta-2>-1.$  We need to study two cases:\\
Case \#1: If $(2-\alpha)q+\beta-2<0$: as $r<t$ and by the change of variable $w=t+r$, we have
\begin{eqnarray*}
(I)&\leq&C_\beta (1-\beta)^{1-q}\int_{0}^{\infty}r^{(2-\alpha)q-2}\int_{r}^{\infty}|u^{(2)}(x,t+r)|^{q}dt\,dr\\
&\leq& C_\beta(1-\beta)^{1-q}\int_{0}^{\infty}r^{[(2-\alpha)q-1]-1}\int_{2r}^{\infty}|u^{(2)}(x,w)|^{q}dw\,dr\\
&\leq&C_\beta (1-\beta)^{1-q}\int_{0}^{\infty}r^{[(2-\alpha)q-1]-1}\int_{r}^{\infty}|u^{(2)}(x,w)|^{q}dw\,dr.
\end{eqnarray*}
Then by Hardy's inequality (\ref{hardy2}) as $(2-\alpha)q-1>0$
\begin{eqnarray*}
(I)&\leq&C_\beta (1-\beta)^{1-q} \frac{1}{(2-\alpha)q-1}\int_{0}^{\infty}\big(w^{2-\alpha}|u^{(2)}(x,w)|\big)^{q}\frac{dw}{w}.
\end{eqnarray*}
Case \#2: If $(2-\alpha)q+\beta-2\geq 0$: By the change of variable $w=t+r$, we get
\begin{eqnarray*}
(I)
&\leq&C_\beta (1-\beta)^{1-q}\int_{0}^{\infty}r^{-\beta} \int_{r}^{\infty}(t+r)^{(2-\alpha)q+\beta-2}|u^{(2)}(x,t+r)|^{q}dt\,dr\\
&=& C_\beta(1-\beta)^{1-q}\int_{0}^{\infty}r^{-\beta} \int_{2r}^{\infty}w^{(2-\alpha)q+\beta-2}|u^{(2)}(x,w)|^{q} dw\,dr\\
&\leq& C_\beta(1-\beta)^{1-q}\int_{0}^{\infty}r^{-\beta} \int_{r}^{\infty}w^{(2-\alpha)q+\beta-2}|u^{(2)}(x,w)|^{q} dw\,dr,\\
\end{eqnarray*}
and by Hardy's inequality (\ref{hardy2}),
\begin{eqnarray*}
(I)&\leq&\frac{C_\beta}{(1-\beta)^{q}}\int_{0}^{\infty}\big(w^{2-\alpha}|u^{(2)}(x,w)|\big)^{q}\frac{dw}{w}.
\end{eqnarray*}
Therefore, in both cases we have
\begin{eqnarray*}
(I)&\leq&C_\beta \int_{0}^{\infty}\big(w^{2-\alpha}|u^{(2)}(x,w)|\big)^{q}\frac{dw}{w}.
\end{eqnarray*}

To estimate $(II)$ observe  $r^{-\beta}<t^{-\beta}$, for $r>t$ and $\beta>0$, then we use the same argument used before to estimate $(I)$ case \#1, doing the change of variable $w=t+r$, and using we get  Hardy's inequality (\ref{hardy2}) to get
\begin{eqnarray*}
(II)&\leq&\frac{C_\beta}{1-\alpha}\int_{0}^{\infty}\big(w^{2-\alpha}|u^{(2)}(x,w)|\big)^{q}\frac{dw}{w}.
\end{eqnarray*}
Finally,
$$\|\big(\int_{0}^{\infty}\big(t^{1-(\alpha-\beta)}|\frac{\partial
}{\partial t}P_{t}(D^{\beta}f)|\big)^{q}\frac{dt}{t}\big)^{1/q}\|_{p,\gamma} \hspace{4.0cm}$$
\begin{eqnarray}\label{est11}
\nonumber &&\hspace{2.0cm}  \leq C \|
(\int_{0}^{\infty}\big(t^{1-(\alpha-\beta)}\int_{0}^{+\infty}r^{-\beta}|u^{(2)}(\cdot,t+r)|dr\big)^{q}\frac{dt}{t})^{1/q}\|_{p,\gamma}\\
&& \hspace{2.0cm} \leq C \|\big(\int_{0}^{\infty}\big(w^{2-\alpha}|u^{(2)}(\cdot,w)|\big)^{q}\frac{dw}{w}\big)^{1/q}\|_{p,\gamma}<\infty,
\end{eqnarray}
as $f\in F_{p,q}^{\alpha}(\gamma_{d}).$ By the previous estimate and (\ref{est5})
$$
\|D^{\beta}f\|_{F_{p,q}^{\alpha-\beta}}\leq C \|f\|_{F_{p,q}^{\alpha}}.$$
\ep
\\

In following theorem we will study the boundedness of the Bessel fractional derivative on Triebel-Lizorkin spaces, for the case $0<\beta<\alpha<1$,
\begin{teo}\label{DerBessel<1}
Let $0<\beta<\alpha<1$, $1\leq p, q <\infty$ then
 ${\mathcal D}^{\beta}$ is
bounded from $F_{p,q}^{\alpha}(\gamma_{d})$ into
$F_{p,q}^{\alpha-\beta}(\gamma_{d})$.
\end{teo}

\dem

Let $f\in L^{p}(\gamma_{d})$, using the Fundamental Theorem of Calculus we can write,
 \begin{eqnarray*}
\nonumber |\mathcal{D}^{\beta}f(x)| &\leq&\displaystyle\frac{1}{c_{\beta}}\int_{0}^{+\infty}s^{-\beta-1}|e^{-s}P_{s}f(x)-f(x)|ds\\
&\leq&\displaystyle\frac{1}{c_{\beta}}\int_{0}^{+\infty}s^{-\beta-1}e^{-s}|P_{s}f(x)-f(x)|ds
+\displaystyle\frac{1}{c_{\beta}}\int_{0}^{+\infty}s^{-\beta-1}|e^{-s}-1||f(x)|ds\\
&\leq&\displaystyle\frac{1}{c_{\beta}}\int_{0}^{+\infty}s^{-\beta-1}\int_{0}^{s}|\frac{\partial
}{\partial r}P_{r}f(x)|dr\,ds
+\displaystyle\frac{1}{c_{\beta}} |f(x)|\int_{0}^{+\infty}s^{-\beta-1}\int_{0}^{s}e^{-r}dr\,ds.
\end{eqnarray*}

Now, using Hardy's inequality  (\ref{hardy1}), with $p=1$, in both integrals,  we have
 \begin{eqnarray*}
\nonumber |\mathcal{D}^{\beta}f(x)|
&\leq&\displaystyle\frac{1}{c_{\beta}}\int_{0}^{+\infty}s^{-\beta-1}\int_{0}^{s}|\frac{\partial
}{\partial r}P_{r}f(x)|dr\,ds
+\displaystyle\frac{1}{c_{\beta}}|f(x)|\int_{0}^{+\infty}s^{-\beta-1}\int_{0}^{s}e^{-r}dr \,ds\\
&\leq& \displaystyle\frac{1}{\beta c_{\beta}}\int_{0}^{+\infty}r^{1-\beta}|\frac{\partial
}{\partial r}P_{r}f(x)|\frac{dr}{r}
+\displaystyle\frac{1}{\beta c_{\beta}}\Gamma(1-\beta)|f(x)|.
\end{eqnarray*}

Thus,
$$ |\mathcal{D}^{\beta}f(x)| \leq   \displaystyle\frac{1}{\beta c_{\beta}}\int_{0}^{+\infty}r^{1-\beta}|u^{(1)}(x,r)|\frac{dr}{r}
+\displaystyle\frac{1}{\beta c_{\beta}}\Gamma(1-\beta)|f(x)|.$$

Therefore, if  $f\in F_{p,q}^{\alpha}(\gamma_{d})$, we get
 \begin{eqnarray}\label{contderbes}
\nonumber\|\mathcal{D}^{\beta}f\|_{p,\gamma}&\leq&
\frac{1}{\beta c_{\beta}}\|\displaystyle\int_{0}^{+\infty} r^{1-\beta}|u^{(1)}(\cdot,r)|\frac{dr}{r}\|_{p,\gamma}
+\displaystyle\frac{1}{\beta c_{\beta}}\Gamma(1-\beta)\|f\|_{p,\gamma}\\
&\leq&C_{\beta}\|f\|_{F_{p,1}^{\beta}}\leq C'_{\beta}\|f\|_{F_{p,q}^{\alpha}},
\end{eqnarray}
since  $ F_{p,q}^{\alpha}(\gamma_{d})\subset
F_{p,1}^{\beta}(\gamma_{d}) $, as $\alpha> \beta,$ and $q\geq1$. \\

On the other hand, using a similar argument as above (the Fundamental Theorem of Calculus and Hardy's inequality  (\ref{hardy1}) with $p=1$),  we have,
 \begin{eqnarray*}
|\frac{\partial}{\partial t} P_{t}(\mathcal{D}^{\beta}f )(x)|&\leq&\frac 1{c_\beta}\int_0^{\infty}s^{-\beta-1} |e^{-s}\frac{\partial}{\partial t} P_{t+s}f(x)  -\frac{\partial}{\partial t} P_tf(x)| ds\\
&\leq&\frac 1{c_\beta}\int_0^{\infty}s^{-\beta-1}e^{-s}| \frac{\partial}{\partial t} P_{t+s}f(x)  -\frac{\partial}{\partial t} P_tf(x)| ds\\
&& \quad \quad +\frac 1{c_\beta}\int_0^{\infty}s^{-\beta-1}|e^{-s}  -1| |\frac{\partial}{\partial t} P_{t}f(x)| ds\\
&\leq&\frac 1{c_\beta}\int_0^{\infty}s^{-\beta-1} \int_0^{s} 
|u^{(2)}(x,t+r)| dr \,ds\\
&& \quad  \quad + \frac 1{c_\beta} |u^{(1)}(x,t)| \int_0^{\infty}s^{-\beta-1}\int_{0}^{s}e^{-r}dr\,ds,\\
&\leq& \frac 1{\beta c_\beta}\int_0^{\infty}r^{-\beta}  
|u^{(2)}(x,t+r)| dr +\frac {\Gamma(1-\beta)}{\beta c_\beta}|u^{(1)}(x,t)|.
\end{eqnarray*}
Therefore, 
$$
|\frac{\partial}{\partial t} P_{t}(\mathcal{D}^{\beta}f (x))| \leq \frac 1{\beta c_\beta}\int_0^{\infty}r^{-\beta}  
|u^{(2)}(x,t+r)| dr +\frac {\Gamma(1-\beta)}{\beta c_\beta}|u^{(1)}(x,t)|,$$
and then, we have
$$\|\big(\int_{0}^{\infty}\big(t^{1-(\alpha-\beta)}|\frac{\partial }{\partial t}P_{t}(\mathcal{D}^{\beta}f)|\big)^{q}\frac{dt}{t}\big)^{1/q}\|_{p,\gamma}\hspace{5.0cm}$$
\begin{eqnarray*}
&\leq&\frac {C}{\beta c_\beta}\|\big(\int_{0}^{\infty}\big(t^{1-(\alpha-\beta)}\int_0^{\infty}r^{-\beta}  
|u^{(2)}(\cdot,t+r)| dr \big)^{q}\frac{dt}{t}\big)^{1/q}\|_{p,\gamma}\\
&&\hspace{4.0cm} +\frac {C}{\beta c_\beta}\Gamma(1-\beta)\|\big(\int_{0}^{\infty}\big(t^{1-(\alpha-\beta)}|u^{(1)}(\cdot,t)|\big)^{q}\frac{dt}{t}\big)^{1/q}\|_{p,\gamma}.
\end{eqnarray*}
Since the exponential factor has been remove, the first term can be estimate as in the proof of Theorem 2.3, estimates (\ref{est10}) and (\ref{est11}),
\begin{eqnarray*}
\|\big(\int_{0}^{\infty}\big(t^{1-(\alpha-\beta)}\int_0^{\infty}r^{-\beta}  
|u^{(2)}(\cdot,t+r)| dr \big)^{q}\frac{dt}{t}\big)^{1/q}\|_{p,\gamma}&\leq&C\|\int_{0}^{\infty}\big(t^{2-\alpha}|\frac{\partial^{2}
}{\partial
t^{2}}P_{t}f|\big)^{q}\frac{dt}{t}\big)^{1/q}\|_{p,\gamma}\\
\end{eqnarray*}
which is finite as $f\in F_{p,q}^{\alpha}(\gamma_{d}),$ 
and for the second term, we have
\begin{eqnarray*}
\|\big(\int_{0}^{\infty}\big(t^{1-(\alpha-\beta)}|u^{(1)}(x,t)|\big)^{q}\frac{dt}{t}\big)^{1/q}\|_{p,\gamma}&\leq&C\|f\|_{F_{p,q}^{\alpha-\beta}}\leq C\|f\|_{F_{p,q}^{\alpha}}
\end{eqnarray*}
as  $ F_{p,q}^{\alpha}(\gamma_{d})\subset
F_{p,q}^{\alpha-\beta}(\gamma_{d}) $,  
thus, 
\begin{eqnarray*}
\|\big(\int_{0}^{\infty}\big(t^{1-(\alpha-\beta)}|\frac{\partial }{\partial t}P_{t}(\mathcal{D}^{\beta}f)|\big)^{q}\frac{dt}{t}\big)^{1/q}\|_{p,\gamma}
&\leq C\|f\|_{F_{p,q}^{\alpha}}.
\end{eqnarray*}
Therefore, $\mathcal{D}^{\beta}f\in
F_{p,q}^{\alpha-\beta}(\gamma_{d})$ and moreover,  by previous estimate and (\ref{contderbes})
\begin{eqnarray*}
\|\mathcal{D}^{\beta}f\|_{F_{p,q}^{\alpha-\beta}}
&\leq&C\|f\|_{F_{p,q}^{\alpha}}. 
\end{eqnarray*}
 \ep

To study  the general case for fractional derivatives (removing the condition that the indexes must be less than 1), we need to consider forward differences. For a given function $f$,  the $k$-th order forward difference of $f$ starting at $t$ with increment $s$ is defined as,\\
 $$\Delta_{s}^{k}(f,t)=\displaystyle\sum_{j=0}^{k}\binom{k}{j}(-1)^{j}f(t+(k-j)s).$$
 
 We will need the following technical results
 \begin{lemma}\label{forw-diff} For any positive integer $k$
 \begin{enumerate}
 \item[i)]$\Delta_{s}^{k}(f,t)=\Delta_{s}^{k-1}(\Delta_{s}(f,\cdot),t)=\Delta_{s}(\Delta_{s}^{k-1}(f,\cdot),t)$
 \item[ii)] $\Delta_{s}^{k}(f,t)=\displaystyle\int_{t}^{t+s}\int_{v_{1}}^{v_{1}+s}...
 \int_{v_{k-2}}^{v_{k-2}+s}\int_{v_{k-1}}^{v_{k-1}+s}f^{(k)}(v_{k})dv_{k}dv_{k-1}...dv_{2}dv_{1}$
\item [iii)]For any positive integer $k$,
\begin{equation}\label{difder}
\frac{\partial }{\partial s}(\Delta_s^k (f,t) ) = k \,\Delta_s^{k-1} (f',t+s),
\end{equation}
and for any integer $j>0$,
\begin{equation}\label{difder2}
\frac{\partial^j }{\partial t^j}(\Delta_s^k (f,t) ) =\Delta_s^{k} (f^{(j)},t).
\end{equation} \end{enumerate}
 \end{lemma}
 The proof of this known lemma can also be found in an appendix in \cite{gatur}.
\begin{lemma}\label{hardyineqk}(Hardy's type innequality)
Let $f\geq 0, r>0, p\geq 1$ and $k\in\mathbb{N}$ then
 \begin{eqnarray*}
&&\big(\int_{0}^{+\infty}\big(\int_{0}^{x}...\int_{0}^{x}f(r_{1},...,r_{k})\,dr_{1}...dr_{k})\big)^{p}x^{-r-1}dx\big)^{1/p}\\
&&\leq\int_{0}^{1}...\int_{0}^{1}\big(\int_{0}^{+\infty}(x^{k}f(xv_{1},...,xv_{k}))^{p}x^{-r-1}dx\big)^{1/p}dv_{1}...dv_{k}
 \end{eqnarray*}
\end{lemma}
\dem
Taking $r_{1}=xv_{1},...,r_{k}=xv_{k}$, we get
 \begin{eqnarray*}
 &&\big(\int_{0}^{+\infty}\big(\int_{0}^{x}...\int_{0}^{x}f(r_{1},...,r_{k})dr_{1}...dr_{k})\big)^{p}x^{-r-1}dx\big)^{1/p}\\
 &=&\big(\int_{0}^{+\infty}\big(\int_{0}^{1}...\int_{0}^{1}f(xv_{1},...,xv_{k})x^{k}dv_{1}...dv_{k})\big)^{p}x^{-r-1}dx\big)^{1/p}\\
&=&\big(\int_{0}^{+\infty}\big(\int_{(0,1)^{k}}f(xv_{1},...,xv_{k})x^{k}dv_{1}...dv_{k})\big)^{p}x^{-r-1}dx\big)^{1/p}
\end{eqnarray*}
Now, consider the spaces $L^{p}((0,+\infty),x^{-r-1})$ and $L^{p}((0,1)^{k})$; then by Minkowski's inequality,
 \begin{eqnarray*}
 &&\big(\int_{0}^{+\infty}\big(\int_{(0,1)^{k}}f(xv_{1},...,xv_{k})x^{k}dv_{1}...dv_{k})\big)^{p}x^{-r-1}dx\big)^{1/p}\\
&\leq&\int_{(0,1)^{k}}\big(\int_{0}^{+\infty}\big(f(xv_{1},...,xv_{k})x^{k}\big)^{p}x^{-r-1}dx\big)^{1/p}dv_{1}...dv_{k}\\
&=&\int_{0}^{1}...\int_{0}^{1}\big(\int_{0}^{+\infty}\big(x^{k}f(xv_{1},...,xv_{k})\big)^{p}x^{-r-1}dx\big)^{1/p}dv_{1}...dv_{k}.
\end{eqnarray*}
Thus,
 \begin{eqnarray*}
 &&\big(\int_{0}^{+\infty}\big(\int_{0}^{x}...\int_{0}^{x}f(r_{1},...,r_{k})dr_{1}...dr_{k})\big)^{p}x^{-r-1}dx\big)^{1/p}\\
&\leq&\int_{0}^{1}...\int_{0}^{1}\big(\int_{0}^{+\infty}\big(x^{k}f(xv_{1},...,xv_{k})\big)^{p}x^{-r-1}dx\big)^{1/p}dv_{1}...dv_{k}
\end{eqnarray*}
\quad \quad\quad         \ep

 \begin{lemma}\label{forw-diff-TL} For any positive integer $k$,
  \begin{eqnarray*}
  \Delta_{s}^{k}(f,t)&=&\int_{t}^{t+s}\int_{v_{1}}^{v_{1}+s}...
 \int_{v_{k-1}}^{v_{k-1}+s}f^{(k)}(v_{k})dv_{k}...dv_{2}dv_{1}\\
&=&\int_{0}^{s}...
\int_{0}^{s}f^{(k)}(t+v_{1}+...+v_{k})dv_{k}...dv_{1}
\end{eqnarray*}
 \end{lemma}
  \begin{lemma}\label{fdsk} Let $t\geq 0, \beta>0$ and let $k$ be the smallest integer greater than $\beta$ and let $f$ differentiable up to order $k$, then
  \begin{eqnarray*}
  \int_{0}^{+\infty}s^{-\beta-1}|\Delta_{s}^{k}(f,t)|ds\leq C_{\beta,k} \int_{0}^{+\infty}w^{k-\beta-1}|f^{(k)}(t+w)|dw
  \end{eqnarray*}
  where $C_{\beta,k}=\displaystyle\int_{0}^{1}...\int_{0}^{1}(v_{1}+...+v_{k})^{\beta-k}dv_{1}...dv_{k}$
   \end{lemma}
   \dem
From Lemmas \ref{hardyineqk} and \ref{forw-diff-TL} we have,
    \begin{eqnarray*}
    \int_{0}^{+\infty}s^{-\beta-1} |\Delta_{s}^{k}(f,t)|ds&\leq& \int_{0}^{+\infty}s^{-\beta-1}\int_{0}^{s}...
 \int_{0}^{s}|f^{(k)}(t+v_{1}+...+v_{k})|dv_{1}...dv_{k}ds\\
 &\leq&\int_{0}^{1}...\int_{0}^{1}\big(\int_{0}^{+\infty}(s^{k}|f^{(k)}(t+s(v_{1}+...+v_{k})|)s^{-\beta-1}ds\big)dv_{1}...dv_{k}\\
 &=&\int_{0}^{1}...\int_{0}^{1}\big(\int_{0}^{+\infty}(s^{k-\beta-1}|f^{(k)}(t+s(v_{1}+...+v_{k})|)ds\big)dv_{1}...dv_{k}
  \end{eqnarray*}
  taking $r=s(v_{1}+...+v_{k})$ then $dr=(v_{1}+...+v_{k})ds$,
  \begin{eqnarray*}
   \int_{0}^{+\infty}s^{k-\beta-1}|f^{(k)}(t+s(v_{1}+...+v_{k})|ds&=&\int_{0}^{+\infty}r^{k-\beta}|f^{(k)}(t+r)|\frac{dr}{r}(v_{1}+...+v_{k})^{\beta-k}\\
   &=&\int_{0}^{+\infty}r^{k-\beta-1}|f^{(k)}(t+r)|dr(v_{1}+...+v_{k})^{\beta-k}.
   \end{eqnarray*}
 Therefore,
    \begin{eqnarray*}
   \int_{0}^{+\infty}s^{-\beta-1} |\Delta_{s}^{k}(f,t)|ds \quad \quad \quad \quad \quad \quad \quad \quad \quad \quad \quad \quad \quad \quad \quad \quad
  \end{eqnarray*}
  \begin{eqnarray*}
&\leq&\int_{0}^{1}...\int_{0}^{1}\big(\int_{0}^{+\infty}r^{k-\beta-1}|f^{(k)}(t+r)|dr(v_{1}+...+v_{k})^{\beta-k}\big)dv_{1}...dv_{k}\\
     &=&\big(\int_{0}^{+\infty}r^{k-\beta-1}|f^{(k)}(t+r)|dr\big)\int_{0}^{1}...\int_{0}^{1}(v_{1}+...+v_{k})^{\beta-k}dv_{1}...dv_{k}\\
   &=&C_{\beta,k}\big(\int_{0}^{+\infty}r^{k-\beta-1}|f^{(k)}(t+r)|dr,
 \end{eqnarray*}
    where $C_{\beta,k}=\displaystyle\int_{0}^{1}...\int_{0}^{1}(v_{1}+...+v_{k})^{\beta-k}dv_{1}...dv_{k}<\infty.$ \ep

\begin{obs}\label{Binomial}
Using the Binomial Theorem and the semigroup property of $\{P_s\}$, we have 
 \begin{eqnarray}\label{powerrep}
\nonumber ( P_s-I )^k f(x) &=& \sum_{j=0}^k {k \choose j} P_s^{k-j} (-I)^j f(x) = \sum_{j=0}^k {k \choose j} (-1)^jP_s^{k-j} f(x)\\
\nonumber &=&\sum_{j=0}^k {k \choose j} (-1)^jP_{(k-j)s} f(x) =\sum_{j=0}^k {k \choose j} (-1)^j u(x,(k-j)s)\\
&=& \Delta_s^k (u(x, \cdot), 0),
\end{eqnarray}
where  as usual, $u(x,s) = P_s f(x)$and furthermore,
\begin{equation}\label{compo}
P_t  ( P_s-I )^k f(x)  =  \Delta_s^k (u(x, \cdot), t).
\end{equation}

\end{obs}

The following result extends Theorem 2.3 to the general case $0<\beta<\alpha$,
\begin{teo}\label{DerRiesz>1}
Let $0<\beta<\alpha$, $1\leq p, q < \infty$ then $D^{\beta}$ is
bounded from $F_{p,q}^{\alpha}(\gamma_{d})$ into
$F_{p,q}^{\alpha-\beta}(\gamma_{d})$.
\end{teo}
\dem
Let $f\in F_{p,q}^{\alpha}(\gamma_{d})$, using Observation \ref{Binomial} and Lemma \ref{fdsk}
 \begin{eqnarray*}
  |D_{\beta}f(x)|&\leq&\displaystyle\frac{1}{c_{\beta}}\int_{0}^{+\infty}s^{-\beta-1}|(P_{s}-I)^{k}f(x)|ds\\
&=&\displaystyle\frac{1}{c_{\beta}}\int_{0}^{+\infty}s^{-\beta-1}|\Delta_{s}^{k}(u(x,\cdot),0)|ds\\
&\leq&C_{\beta,k}\int_{0}^{+\infty}r^{k-\beta-1} |u^{(k)}(x,r)|dr,
\end{eqnarray*}
then
\begin{eqnarray*}
  \|D_{\beta}f\|_{p,\gamma}&\leq&C_{\beta,k}\|\int_{0}^{+\infty}r^{k-\beta} |u^{(k)}(\cdot,r)|\frac{dr}{r}\|_{p,\gamma} < \infty,
  \end{eqnarray*}
  since $ F_{p,q}^{\alpha}(\gamma_{d})\subset
F_{p,1}^{\beta}(\gamma_{d})$, ($\alpha>\beta$ and $1 \leq q<\infty$).\\

On the other hand, let $n\in\mathbb{N}, n>\alpha$; by Observation \ref{Binomial} and Lemma \ref{fdsk}, we get,
\begin{eqnarray*}
|\frac{\partial^{n}}{\partial t^{n}}P_{t}(D_{\beta}f)(x)|&\leq& \frac{1}{c_{\beta}}\int_{0}^{+\infty}s^{-\beta-1} |\Delta_{s}^{k}(u^{(n)}(x,\cdot),t)|ds\\
&\leq& \frac{1}{c_{\beta}}C_{\beta,k}\int_{0}^{+\infty}r^{k-\beta-1} |u^{(n+k)}(x,t+r)|dr.
 \end{eqnarray*}
 The rest of the proof follows the argument used in Theorem 2.3. We will write the details for the sake of completeness.
 
 \begin{eqnarray*}
&& \big(\int_{0}^{\infty}\big(t^{n-(\alpha-\beta)}|\frac{\partial^{n}}{\partial t^{n}}P_{t}(D_{\beta}f)(x)|\big)^{q}\frac{dt}{t}\big)^{1/q}\\
&\leq& C_{\beta,k} \big(\int_{0}^{\infty}\big(t^{n-(\alpha-\beta)}\int_{0}^{+\infty}r^{k-\beta-1} |u^{(n+k)}(x,t+r)|dr\big)^{q}\frac{dt}{t}\big)^{1/q}\\
&\leq& C_{\beta,k} \big(\int_{0}^{\infty}\big(t^{n-(\alpha-\beta)}\int_{0}^{t}r^{k-\beta-1} |u^{(n+k)}(x,t+r)|dr\big)^{q}\frac{dt}{t}\big)^{1/q}\\
&&\quad \quad \quad + C_{\beta,k} \big(\int_{0}^{\infty}\big(t^{n-(\alpha-\beta)}\int_{t}^{+\infty}r^{k-\beta-1} |u^{(n+k)}(x,t+r)|dr \big)^{q}\frac{dt}{t}\big)^{1/q}\\
&\leq& (I)+(II).
\end{eqnarray*}

Writing $s^{\beta-1} = s^{\frac{\beta -1}{q}+ \frac{\beta -1}{q'}},\, \frac{1}{q}+\frac{1}{q'} =1,$ and using H\"older's inequality in the internal integral and Fubini-Tonelli's theorem,
\begin{eqnarray*}
(I)
&\leq&\frac{(k-\beta)^{1/q}}{k-\beta}\big(\int_{0}^{\infty}\big(t^{(n+k-\alpha)q}t^{\beta-k-1}\big(\int_{0}^{t}|u^{(n+k)}(x,t+r)|^{q}r^{k-\beta-1}dr\big)dt\big)^{1/q}\\
&=&C\big(\int_{0}^{\infty}r^{k-\beta-1}\big(\int_{r}^{+\infty}t^{(n+k-\alpha)q+\beta-k-1}|u^{(n+k)}(x,t+r)|^{q}dt\big)dr\big)^{1/q}\\
\end{eqnarray*}

Now, it is easy to see that $(n+k-\alpha)q+\beta-k-1>-1$. Then we need to study two cases:\\

Case\# 1: If $(n+k-\alpha)q+\beta-k-1<0$ then as $w<t$ and taking the change of variable $v=r+t$, 
\begin{eqnarray*}
(I)&\leq&(k-\beta)^{1/q-1}\big(\int_{0}^{\infty}r^{((n+k-\alpha)q-1)-1}\int_{2r}^{\infty}|u^{(n+k)}(x,v)|^{q}dv\big)dr\big)^{1/q}\\
&\leq&(k-\beta)^{1/q-1}\big(\int_{0}^{\infty}r^{((n+k-\alpha)q-1)-1}\int_{r}^{\infty}|u^{(n+k)}(x,v)|^{q}dv\big)dr\big)^{1/q}.
\end{eqnarray*}
Therefore, by Hardy's inequality as $(n+k-\alpha)q-1>0$
\begin{eqnarray*}
(I)&\leq&(k-\beta)^{1/q-1}\big(\frac{1}{(n+k-\alpha)q-1}\big)^{1/q}\big(\int_{0}^{\infty}\big(s^{n+k-\alpha}|u^{(n+k)}(x,s)|\big)^{q}\frac{ds}{s}\big)^{1/q}
\end{eqnarray*}

Case\# 2: If $(n+k-\alpha)q+\beta-k-1\geq 0$ then, making the change of variable $v=r+t$ we get by Hardy's inequality, as $k-\beta>0$,
\begin{eqnarray*}
(I)&\leq&(k-\beta)^{1/q-1}\big(\int_{0}^{\infty}r^{k-\beta-1}\big(\int_{r}^{+\infty}t^{(n+k-\alpha)q+\beta-k-1}|u^{(n+k)}(x,t+r)|^{q}dt\big)dr\big)^{1/q}\\
&\leq&(k-\beta)^{1/q-1}\big(\int_{0}^{\infty}r^{k-\beta-1}\big(\int_{r}^{+\infty}(r+t)^{(n+k-\alpha)q+\beta-k-1}|u^{(n+k)}(x,t+r)|^{q}dt\big)dr\big)^{1/q}\\
&=&(k-\beta)^{1/q-1}\big(\int_{0}^{\infty}r^{k-\beta-1}\big(\int_{2r}^{+\infty}v^{(n+k-\alpha)q+\beta-k-1}|u^{(n+k)}(x,v)|^{q}dv\big)dr\big)^{1/q}\\
&\leq&(k-\beta)^{1/q-1}\big(\int_{0}^{\infty}r^{k-\beta-1}\big(\int_{r}^{+\infty}s^{(n+k-\alpha)q+\beta-k-1}|u^{(n+k)}(x,s)|^{q}ds\big)dr\big)^{1/q}\\
&\leq&\frac{1}{k-\beta}\big(\int_{0}^{\infty}\big( s^{n+k-\alpha}|u^{(n+k)}(x,s)|\big)^{q}\frac{ds}{s}\big)^{1/q}.
\end{eqnarray*}

Thus in both cases, we have
\begin{eqnarray*}
(I)&\leq&C\big(\int_{0}^{\infty}\big(s^{n+k-\alpha}|u^{(n+k)}(x,s)|\big)^{q}\frac{ds}{s}\big)^{1/q}.
\end{eqnarray*}

As $r^{k-\beta-1}\leq t^{k-\beta-1}$ for $r\geq t$, since $k-\beta-1\leq 0$, and taking the change of variable $v=t+r$,
\begin{eqnarray*}
(II)
&\leq&C( \int_{0}^{\infty}\big(t^{n+k-\alpha-1}\int_{t}^{+\infty}|u^{(n+k)}(x,t+r)|dr\big)^{q}\frac{dt}{t}\big)^{1/q}\\
&=&C (\int_{0}^{\infty}\big(t^{n+k-\alpha-1}\int_{2t}^{+\infty}|u^{(n+k)}(x,s)|ds\big)^{q}\frac{dt}{t}\big)^{1/q}\\
&\leq&C (\int_{0}^{\infty}t^{(n+k-\alpha-1)q-1}\big(\int_{t}^{+\infty}|u^{(n+k)}(x,s)|ds\big)^{q}dt\big)^{1/q},
\end{eqnarray*}
then, by Hardy's inequality
\begin{eqnarray*}
(II)&\leq&\frac{C}{n+k-\alpha-1}\big(\int_{0}^{\infty}\big(s^{(n+k-\alpha)}|u^{(n+k)}(x,s)|\big)^{q}\frac{ds}{s}\big)^{1/q}.
\end{eqnarray*}
Therefore,
$$ \|\big(\int_{0}^{\infty}\big(t^{n-(\alpha-\beta)}|\frac{\partial^{n}}{\partial t^{n}}P_{t}D_{\beta}f|\big)^{q}\frac{dt}{t}\big)^{1/q}\|_{p,\gamma} \hspace{23.0cm}$$
\begin{eqnarray}\label{arg10}
&&\nonumber  \hspace{2.5cm} \leq C \|\big(\int_{0}^{\infty}\big(t^{n-(\alpha-\beta)}\int_{0}^{+\infty}r^{k-\beta-1} |u^{(n+k)}(\cdot,t+r)|dr\big)^{q}\frac{dt}{t}\big)^{1/q}\|_{p,\gamma}\\
&&  \hspace{2.5cm} \leq C \|\big(\int_{0}^{\infty}\big(s^{n+k-\alpha}|u^{(n+k)}(\cdot,s)|\big)^{q}\frac{ds}{s}\big)^{1/q}\|_{p,\gamma}\\
&&\nonumber  \hspace{2.0cm} \leq C\|f\|_{F_{p,q}^{\alpha}}.
\end{eqnarray}
Therefore,
$D^{\beta}f\in F_{p,q}^{\alpha-\beta}(\gamma_{d})$ and moreover,
\begin{eqnarray*}
\|D^{\beta}f\|_{F_{p,q}^{\alpha-\beta}}&\leq&C\|f\|_{F_{p,q}^{\alpha}}.
\end{eqnarray*}
\ep 

Finally, the following result extends Theorem 2.4 to the general case $0<\beta<\alpha$,
\begin{teo}\label{DerBessel>1}
Let $0<\beta<\alpha$, $1< p <\infty$ and $1\leq q<\infty,$  then $\mathcal{D}^{\beta}$ is
bounded from $F_{p,q}^{\alpha}(\gamma_{d})$ into
$F_{p,q}^{\alpha-\beta}(\gamma_{d})$.
\end{teo}
\dem
 Let $f\in F_{p,q}^{\alpha}(\gamma_{d})$ and $v(x,r)=e^{-r}u(x,r)$, using Lemma \ref{fdsk} and Leibnitz's differentiation rule for the product
 \begin{eqnarray*}
  |\mathcal{D}^{\beta}f(x)|&\leq&\displaystyle\frac{1}{c_{\beta}}\int_{0}^{+\infty}s^{-\beta-1}|(e^{-s}P_{s}-I)^{k}f(x)|ds=\displaystyle\frac{1}{c_{\beta}}\int_{0}^{+\infty}s^{-\beta-1}|\Delta_{s}^{k}(v(x,\cdot),0)|ds\\
&\leq&C_{\beta,k}\int_{0}^{+\infty}r^{k-\beta} |v^{(k)}(x,r)|\frac{dr}{r}\leq C_{\beta,k}\big(\sum_{j=0}^{k}\binom{k}{j}\int_{0}^{+\infty}r^{k-\beta} e^{-r}|u^{(k-j)}(x,r)|\frac{dr}{r}\big)\\
&=&C_{\beta,k}\big(\sum_{j=0}^{k-1}\binom{k}{j}\int_{0}^{+\infty}r^{k-\beta} e^{-r}|u^{(k-j)}(x,r)|\frac{dr}{r}\big)+ C_{\beta,k}\int_{0}^{+\infty}r^{k-\beta} e^{-r}|u(x,r)|\frac{dr}{r},
\end{eqnarray*}
then
\begin{eqnarray*}
  \|\mathcal{D}^{\beta}f\|_{p,\gamma}&\leq&C_{\beta,k}\big(\sum_{j=0}^{k-1}\binom{k}{j}\|\int_{0}^{+\infty}r^{k-\beta} |u^{(k-j)}(\cdot,r)|\frac{dr}{r}\|_{p,\gamma}\big)+C_{\beta,k}\|\int_{0}^{+\infty}r^{k-\beta} e^{-r}|u(\cdot,r)|\frac{dr}{r}\|_{p,\gamma}\\
  &\leq&C_{\beta,k}\big(\sum_{j=0}^{k-1}\binom{k}{j}\|\int_{0}^{+\infty}r^{k-\beta} |u^{(k-j)}(\cdot,r)|\frac{dr}{r}\|_{p,\gamma}\big)+C_{\beta,k}\int_{0}^{+\infty}r^{k-\beta} e^{-r}\|u(\cdot,r)\|_{p,\gamma}\frac{dr}{r}\\
  &\leq&C_{\beta,k}\big(\sum_{j=0}^{k-1}\binom{k}{j}\|\int_{0}^{+\infty}r^{k-j-(\beta-j)} |u^{(k-j)}(\cdot,r)|\frac{dr}{r}\|_{p,\gamma}\big)+C_{\beta,k}\|f\|_{p,\gamma}\Gamma(k-\beta)\\
  &\leq&C\| f \|_{F_{p,q}^{\alpha}},
  \end{eqnarray*}
because $F_{p,q}^{\alpha}(\gamma_{d})\subset
F_{p,1}^{\beta-j}(\gamma_{d})$, as $\alpha>\beta\geq\beta-j\geq 0$, for $j = 0,...,k-1$ and $q\geq1$.\\

On the other hand, let $n\in\mathbb{N}, n>\alpha$ and $w(x,t)=e^{-t}u^{(n)}(x,t)$, by Lemma \ref{fdsk} we get\begin{eqnarray*}
|\frac{\partial^{n}}{\partial t^{n}}P_{t}(\mathcal{D}^{\beta}f)(x)|&\leq& \frac{e^{t}}{c_{\beta}}\int_{0}^{+\infty}s^{-\beta-1} |\Delta_{s}^{k}(w(x,\cdot),t)|ds\\
&\leq& e^{t}C_{\beta,k}\int_{0}^{+\infty}s^{k-\beta-1} |w^{(k)}(x,t+s)|ds.
 \end{eqnarray*}
 Now, by Leibnitz rule, $w^{(k)}(x,r)=\displaystyle\sum_{j=0}^{k}\binom{k}{j}(-1)^{j}e^{-r}u^{(k+n-j)}(x,r)$ and 
  then
  $$|w^{(k)}(x,r)|\leq \displaystyle\sum_{j=0}^{k}\binom{k}{j}e^{-r}|u^{(k+n-j)}(x,r)|, $$
for all $r >0$. Thus
 \begin{eqnarray*}
|\frac{\partial^{n}}{\partial t^{n}}P_{t}(\mathcal{D}^{\beta}f)(x)|&\leq& C_{\beta,k}\sum_{j=0}^{k}\binom{k}{j}\int_{0}^{+\infty}s^{k-\beta-1} e^{-s}|u^{(k+n-j)}(x,t+s)|ds.
 \end{eqnarray*}
Therefore,
 \begin{eqnarray*}
&& \big(\int_{0}^{\infty}\big(t^{n-(\alpha-\beta)}|\frac{\partial^{n}}{\partial t^{n}}P_{t}(\mathcal{D}^{\beta}f)(x)|\big)^{q}\frac{dt}{t}\big)^{1/q}\\
&\leq&C_{\beta,k}\sum_{j=0}^{k}\binom{k}{j} \big(\int_{0}^{\infty}\big(t^{n-(\alpha-\beta)}\int_{0}^{+\infty}s^{k-j-(\beta-j)-1} e^{-s}|u^{(k-j+n)}(x,t+s)|ds\big)^{q}\frac{dt}{t}\big)^{1/q}
\end{eqnarray*}
For $0\leq j\leq k-1,$ we  have $\beta-j\geq\beta-(k-1)\geq0$ using the same argument as in the proof of Theorem \ref{DerRiesz>1}, see (\ref{arg10}), using Hardy's inequality we get
\begin{eqnarray*}
\| \big(\int_{0}^{\infty}\big(t^{n-(\alpha-\beta)}\int_{0}^{+\infty}s^{k-j-(\beta-j)-1} e^{-s}|u^{(k+n-j)}(x\cdot,t+s)|ds\big)^{q}\frac{dt}{t}\big)^{1/q}\|_{p,\gamma}&<&C \|f\|_{F_{p,q}^{\alpha}}
\end{eqnarray*}
for $0\leq j\leq k-1.$\\
Unfortunately the remaining case $j=k$  requires and special argument,that uses the following known estimate for the Poisson-Hermite semigroup
\begin{equation}\label{est12}
| \frac{\partial^n}{\partial t^n} P_t f(x)| \leq C T^*f(x) t^{-n},
\end{equation}
where $T^*f$ is the maximal function of the Ornstein-Uhlenbeck semigroup. The proof of this estimate can be found in \cite{piur02}, Lemma 2.1.
$$ \big(\int_{0}^{\infty}\big(t^{n-(\alpha-\beta)}\int_{0}^{+\infty}s^{k-\beta-1} e^{-s}|u^{(n)}(\cdot,t+s)|ds\big)^{q}\frac{dt}{t}\big)^{1/q}\hspace{2.0cm}$$
\begin{eqnarray*}
&\leq& C \big(\int_{0}^{\infty}\big(t^{n-(\alpha-\beta)}\int_{0}^{t}s^{k-\beta-1} e^{-s}|u^{(n)}(x,t+s)|ds\big)^{q}\frac{dt}{t}\big)^{1/q}\\
&&+\, C\big(\int_{0}^{\infty}\big(t^{n-(\alpha-\beta)}\int_{t}^{+\infty} s^{k-\beta-1} e^{-s}|u^{(n)}(x,t+s)|ds\big)^{q}\frac{dt}{t}\big)^{1/q}\\
&=&(I)+(II).
\end{eqnarray*}
We consider first the case  $k\leq\alpha$. The term (I) is estimated as term (I) in the proof of Theorem \ref{DerRiesz>1}.

\begin{eqnarray*}
(I)
&\leq&C\big(\int_{0}^{\infty}\big(v^{n-(\alpha-k)}|u^{(n)}(x,v)|\big)^{q}\frac{dv}{v}\big)^{1/q}.
\end{eqnarray*}
 Since  $\beta\geq k-1,$ making the change of variable $v=t+s$ we get
\begin{eqnarray*}
(II)
&\leq&C\big(\int_{0}^{\infty}t^{(n+k-\alpha-1)q-1}\big(\int_{t}^{+\infty}|u^{(n)}(x,t+s)|ds\big)^{q}dt\big)^{1/q}\\
&=&C\big(\int_{0}^{\infty}t^{(n+k-\alpha-1)q-1}\big(\int_{2t}^{+\infty}|u^{(n)}(x,r)|dr\big)^{q}dt\big)^{1/q}\\
&\leq&C\big(\int_{0}^{\infty}t^{(n+k-\alpha-1)q-1}\big(\int_{t}^{+\infty}|u^{(n)}(x,r)|dr\big)^{q}dt\big)^{1/q}.
\end{eqnarray*}
Therefore, by Hardy's inequality (\ref{hardy2}),
\begin{eqnarray*}
(II)&\leq&\frac{C}{n+k-\alpha-1}\big(\int_{0}^{\infty}\big(r^{n-(\alpha-k)}|u^{(n)}(x,r)|\big)^{q}\frac{dr}{r}\big)^{1/q},
\end{eqnarray*}

Next consider the case $k>\alpha$. In this case, using inequality (\ref{est12}) and Hardy's inequality (\ref{hardy1}), we have
\begin{eqnarray*}
(I)
&\leq&C_{n}|T^{\ast}f(x)|\big(\int_{0}^{\infty}t^{-(\alpha-\beta)q-1}\big(\int_{0}^{t}s^{k-\beta-1} e^{-s}ds\big)^{q}dt\big)^{1/q}\\
&\leq&C_{n}|T^{\ast}f(x)|\frac{1}{\alpha-\beta}\big(\int_{0}^{\infty}s^{(k-\alpha)q-1} e^{-sq}ds\big)^{1/q}\\
&=&C_{n}|T^{\ast}f(x)|\frac{1}{(\alpha-\beta)q^{k-\alpha}}\big(\Gamma((k-\alpha)q)^{1/q}.
\end{eqnarray*}
On the other hand,
\begin{eqnarray*}
(II)
&\leq& \big(\int_{0}^{1}t^{(n+k-\alpha-1)q-1}\big(\int_{t}^{+\infty}e^{-s}|u^{(n)}(x,t+s)|ds\big)^{q}dt\big)^{1/q}\\
&&+\big(\int_{1}^{\infty}t^{(n+k-\alpha-1)q-1}\big(\int_{t}^{+\infty}e^{-s}|u^{(n)}(x,t+s)|ds\big)^{q}dt\big)^{1/q}\\
&=&(III)+(IV).
\end{eqnarray*}
By the usual argument using the change of variable $v=t+s$ and Hardy's inequality (\ref{hardy2}), we get
\begin{eqnarray*}
(III)
&\leq&\big(\int_{0}^{1}t^{(n-1)q-1}\big(\int_{t}^{+\infty}|u^{(n)}(x,t+s)|ds\big)^{q}dt\big)^{1/q}\\
&\leq&\big(\int_{0}^{\infty}t^{(n-1)q-1}\big(\int_{t}^{+\infty}|u^{(n)}(x,t+s)|ds\big)^{q}dt\big)^{1/q}\\
&=&\big(\int_{0}^{\infty}t^{(n-1)q-1}\big(\int_{2t}^{+\infty}|u^{(n)}(x,r)|dr\big)^{q}dt\big)^{1/q}\\
&\leq&\big(\int_{0}^{\infty}t^{(n-1)q-1}\big(\int_{t}^{+\infty}|u^{(n)}(x,r)|dr\big)^{q}dt\big)^{1/q}\\
&\leq&\frac{1}{n-1}\big(\int_{0}^{\infty}\big(r^{n}|u^{(n)}(x,r)|\big)^{q}\frac{dr}{r}\big)^{1/q}.
\end{eqnarray*}
Finally using again inequality (\ref{est12}), we get
\begin{eqnarray*}
(IV)
&\leq&\big(\int_{1}^{\infty}t^{(n+k-\alpha-1)q-1}\big(\int_{t}^{+\infty}e^{-s}C_{n}|T^{\ast}f(x)|t^{-n}ds\big)^{q}dt\big)^{1/q}\\
&=&C_{n}|T^{\ast}f(x)|\big(\int_{1}^{\infty}t^{(k-\alpha-1)q-1}e^{-tq}dt\big)^{1/q}\\
&\leq&C_{n}|T^{\ast}f(x)|\big(\int_{1}^{\infty}t^{(k-\alpha-1)q-1}dt\big)^{1/q}=C_{n}|T^{\ast}f(x)|\big(\frac{1}{(\alpha+1-k)q}\big)^{1/q}.
\end{eqnarray*}

Hence, in both cases, we get that
\begin{eqnarray*}
\|\big(\int_{0}^{\infty}\big(t^{n-(\alpha-\beta)}|\frac{\partial^{n}}{\partial t^{n}}P_{t}(\mathcal{D}^{\beta}f)|\big)^{q}\frac{dt}{t}\big)^{1/q}\|_{p,\gamma}&<&\infty,
\end{eqnarray*}
as $f\in F_{p,q}^{\alpha}(\gamma_{d})$. Therefore,
$\mathcal{D}^{\beta}f\in F_{p,q}^{\alpha-\beta}(\gamma_{d})$ and moreover,
\begin{eqnarray*}
\|\mathcal{D}^{\beta}f\|_{F_{p,q}^{\alpha-\beta}}&\leq&C\|f\|_{F_{p,q}^{\alpha}}.
\end{eqnarray*}
  \ep\\


{\bf Observation}
Let us observe that  if instead of considering the {\em Ornstein-Uhlenbeck operator} (\ref{OUop}) and the {\em Poisson-Hermite semigroup} (\ref{PoissonH}) we consider the {\em  Laguerre differential  operator} in $\mathbb R^d_{+}$,
\begin{equation}
\mathcal{L}^{\alpha} = \sum^{d}_{i=1} \bigg[ x_i \frac{\partial^2}{\partial x^2_i}
+ (\alpha_i + 1 - x_i )  \frac{\partial}{\partial x_i} \bigg],
\end{equation}
and the corresponding {\em Poisson-Laguerre semigroup}, or if we consider the {\em Jacobi differential  operator} in $(-1,1)^d$,
\begin{equation}
\mathcal{L}^{\alpha,\beta} = - \sum^{d}_{i=1} \bigg[ (1-x_i^2)\frac{\partial^2}{\partial x_i^2}
+  (\beta_i -\alpha_i-\left(\alpha_i +\beta_i +2\right)x_i) \frac{\partial}{\partial x_i} \bigg],
\end{equation}
and the corresponding {\em Poisson-Jacobi semigroup} (for details we refer to \cite{ur2}), the arguments are completely analogous.  That is to say, we can defined in analogous manner {\em Laguerre-Triebel-Lizorkin spaces, } and {\em Jacobi-Triebel-Lizorkin spaces} then prove that the corresponding notions of Fractional Integrals and Fractional Derivatives  behave similarly.
In order to see this it is more convenient to use the representation (\ref{PoissonH}) of $P_t$ in terms of the one-sided stable measure $\mu^{(1/2)}_t(ds)$, see \cite{piur02}.


\end{document}